\documentclass[12pt]{amsart}
\usepackage{amsmath}
\usepackage{amsthm}
\usepackage{eucal}
\usepackage{times}
\usepackage{euler}
\usepackage{amssymb}
\usepackage{pst-node}

\author{Piotr \'Sniady}
\thanks{Research supported by State Committee for Scientific Research (KBN)
grant \mbox{2 P03A 007 23}}
\address{Institute of Mathematics,
University of Wroclaw, pl.\ Grunwaldzki 2/4, 50-384 Wroclaw, Poland}
\email{Piotr.Sniady@math.uni.wroc.pl}

\title{Multinomial identities arising from free probability theory}

\theoremstyle{plain}
\newtheorem{lemma}{Lemma}
\newtheorem{theorem}[lemma]{Theorem}
\newtheorem{proposition}[lemma]{Proposition}
\newtheorem{corollary}[lemma]{Corollary}

\theoremstyle{definition}

\theoremstyle{remark}
\newtheorem*{remark}{Remark}
\newtheorem*{example}{Example}

\newcommand{\A}{{\mathcal{A}}}
\newcommand{\E}{{\mathbb{E}}}

\newcommand{\Ha}{{\mathcal{H}}}
\newcommand{\Ka}{{\mathcal{K}}}
\newcommand{\C}{{\mathbb{C}}}
\newcommand{\R}{{\mathbb{R}}}
\newcommand{\Z}{{\mathbb{Z}}}

\newcommand{\N}{{\mathbb{N}}}
\newcommand{\gwia}{^{\star}}

\newcommand{\x}{{\mathbf{x}}}
\newcommand{\D}{{\mathcal{B}}}
\newcommand{\dd}{{\ \mathrm{d}}}
\newcommand{\kk}{{\kappa}}

\DeclareMathOperator{\Tr}{Tr}

\DeclareMathOperator{\vol}{vol}
\DeclareMathOperator{\NC}{NC}

\begin{document}

\bibliographystyle{alpha}

\begin{abstract}
We prove a family of new identities fulfilled by multinomial coefficients,
which were conjectured by Dykema and Haagerup. Our method bases on a study of
the, so--called, triangular operator $T$ by the means of the free probability
theory.
\end{abstract}

\maketitle

\section{Introduction}
\label{sec:intro}
\subsection{Overview}
In order to answer some questions in the theory of operator algebras Dykema and
Haagerup started investigation of the, so--called, triangular operator $T$
\cite{DykemaHaagerup2001}. Currently there are many different descriptions
of this operator: in terms of random matrices, in terms of free probability
theory and a purely combinatorial one (and we will recall them in the following).

Dykema and Haagerup conjectured that the moments of this operator fulfill
\begin{equation}
\phi\big[ \big(T^k
(T\gwia)^k\big)^n\big]=\frac{n^{nk}}{(nk+1)!}
\label{eq:hipoteza}
\end{equation}
for any $k,n\in\N$.
By using the combinatorial description of $T$ Dykema and Yan \cite{DykemaYan}
showed that this
conjecture would imply an infinite family of identities for multinomial
coefficients. For example, for $n=2$ the conjecture is equivalent to a
well--known identity (in the following we will be taking sums over nonnegative
integers):
\begin{equation}
2^{2k}=\sum_{p+q=k} \binom{2p}{p} \binom{2q}{q},
\label{eq:tozsamosc2}
\end{equation}
while for $n=3$ is equivalent to the following one, not known before:
\begin{multline}
\label{eq:tozsamosc3}
3^{3k} =\sum_{p+q=k} \binom{3p}{p,p,p} \binom{3q}{q,q,q} + \\
 + 3 \sum_{\substack{p+q+r=k-1\\ r'+q'=r+q+1 \\ p''+r''=p+r+1}}
\binom{2p+p''}{p,p,p''} \binom{2q+q'}{q,q,q'} \binom{r+r'+r''}{r,r',r''}.
\end{multline}
The complication of the formula grows superexponentially with $n$ and already for $n=4$
it becomes very complicated:
\begin{multline}
\label{eq:tozsamosc4}
4^{4k} =\sum_{p+q=k} \binom{4p}{p,p,p,p} \binom{4q}{q,q,q,q} + \\
 \shoveleft{+8 \!\!
 \sum_{\substack{p+q+r=k-1 \\  p'+q'=p+q+1 \\ p''+q''=p+q+1
\\ q'''+r'''=q+r+1}}
\!\! \binom{2p+p'+p''}{p,p,p',p''}
\binom{q+q'+q''+q'''}{q,q',q'',q'''} \binom{3r+r'''}{r,r,r,r'''}+} \\
\shoveleft{ + 4 \!\!
\sum_{\substack{p+q'+r'=k-1 \\ p+q''+r''=k-1 \\
p'''+q'''=p+q'+1 \\ p''''+q''''=p+q''+1}}
\!\!
\binom{2p+p'''+p''''}{p,p,p''',p''''}
\binom{q'+q'''}{q',q'''} \binom{q''+q''''}{q'',q''''} \times} \\
\shoveright{ \times \binom{2r''}{r'',r''}  \binom{2r'}{r',r'}
\binom{q'+q''+q'''+q''''+2r'+2r''+2}{q'+q'''+2r'+1,q''+q''''+2r''+1}+}\\
\shoveleft{+8 \!\!
\sum_{\substack{p+q+r+s=k-2 \\ q'+r'=q+r+s+2 \\ p''+r''=p+q+r+2}}
\!\! \binom{2p}{p,p} \binom{q+q'}{q,q'} \binom{r+r''}{r,r''} \binom{2s}{s,s}
\times}\\
\shoveright{ \times \binom{3p+p''+2q+q'+2}{2p+q+q'+1,p+q+1,p''}
\binom{2r+r'+r''+3s+2}{r+r''+2s+1,r+s+1,r'}.}
\end{multline}

The main result of this article is the proof of the conjecture
(\ref{eq:hipoteza}). Since this conjecture was originally formulated in the
language of the free probability theory, it should not be a great surprise
that also its proof is formulated in this language. However, in order to make
this article as comprehensible as possible for the wide audience, we will use
only combinatorial methods of free probability and include all necessary
notions.

This article is organized as follows.
In the remaining part of Sect.\ \ref{sec:intro} we present briefly the operator
$T$ from the random matrix and operator theoretic point of view. In Sect.\
\ref{sec:free} we present some rudimentary concepts from the probability
theory, in particular the notion of the generalized circular element. In Sect.\
\ref{sec:main} we define the triangular operator $T$ as a certain
generalized circular elements and we prove the conjecture \eqref{eq:hipoteza}.
In Sect.\ \ref{sec:kombinatoryczna} we 
show how to express moments of $T$ in
terms of multinomial coefficients and hence to get Eq.\ 
\eqref{eq:tozsamosc2}--\eqref{eq:tozsamosc4} (and an infinite number of other identities).
Finally, Sect.\ \ref{sec:techniczne} is devoted to some technical proofs.

\subsection{Triangular operator $T$---motivations}
This section is independent of the rest of this article and can be
skipped by a reader interested only in the combinatorics. We hope, however,
that a broader view presented here might be interesting even for mathematicians
not involved in operator theory.

Von Neumann algebras are algebras of bounded operators on a Hilbert space
$\Ha$, containing identity operator and $A\gwia$ when they contain $A$ and
closed in strong--operator topology \cite{KadisonRingrose1,KadisonRingrose2}.
Especially interesting are type $II_1$ factors $\A$, i.e.\ von Neumann algebras
equipped with a unique normalized faithful tracial state
$\phi:\A\rightarrow\C$.

There are many important open questions concerning type $II_1$ factors, for
example the famous invariant subspace conjecture asking if for every
$x\in\A$ (which is not a multiple of identity operator) there exists a
closed invariant subspace $\Ka\subset\Ha$ which is nontrivial ($\Ka\neq\{0\}$
and $\Ka\neq\Ha$) and such that the orthogonal projection
$\pi_\Ka:\Ha\rightarrow\Ka$ fulfills $\pi_\Ka\in\A$. There are many partial
results concerning this question and the most recent one is due to Haagerup
\cite{Haagerup2001}; he shows that such an invariant subspace exists if
the operator $x$ can be appoximated in a certain way by finite--dimensional
matrices (or, strictly speaking, if $\A$ is embeddable into $R^\omega$, the
ultrapower of the hyperfinite factor) and if the eigenvalues of $x$ are not all
equal (or, strictly speaking, if the Brown measure of $x$ is not supported in a
single point). This result restricts strongly the class of possible
counterexamples for the invariant subspace conjecture and suggests us to study
quasinilpotent operators, for which the assumptions of Haagerup's theorem are
not fulfilled.

The triangular operator $T$ of Dykema and Haagerup arose as a natural
candidate for such a counterexample. The distribution of $T$ was defined
originally \cite{DykemaHaagerup2001} as the limit of distribution of random
matrices $T_N$:
\begin{equation}
\phi(T^{s_1} \cdots T^{s_n})=\lim_{N\rightarrow\infty}
\frac{1}{N} \E \Tr T_N^{s_1} \cdots T_N^{s_n}
\label{eq:oryginalnadefinicja}
\end{equation}
for any $n\in\N$ and $s_1,\dots,s_n\in\{1,\star\}$, where
$$T_N = \left[ \begin{array}{ccccc} t_{1,1} &
t_{1,2} &   \cdots & t_{1,n-1} & t_{1,n} \\ 0 & t_{2,2} & \cdots & t_{2,n-1} &
t_{2,n} \\ \vdots&                 & \ddots & \vdots   & \vdots        \\ & & &
t_{n-1,n-1} & t_{n-1,n} \\ 0&       & \cdots &             0               &
t_{n,n} \end{array} \right]$$
and $(t_{i,j})_{1\leq i\leq j\leq N}$ are independent centered
Gaussian random variables with variance $\frac{1}{N}$. 

When this article was nearly finished Dykema and Haagerup announced that they had proved 
existence of nontrivial hyperinvariant subspaces of the operator $T$. Their proof 
uses Theorem \ref{theo:formulanapochodna} from this article.



\section{Operator--valued free probability theory}
\label{sec:free}
\subsection{Free probability}
Free probability theory was initiated by Voiculescu in order to answer some old
questions in the theory of operator algebras
\cite{VoiculescuDykemaNica,VoiculescuPlenar,VoiculescuLectures}, but it soon
evolved into an exciting self--standing theory with many links to other fields,
to mention only the theory of random matrices \cite{Voiculescu1991}, theory of
representations of groups of permutations $S_n$ \cite{Biane1998} and theoretical
physics \cite{SpeicherNeu}. Free probability has also its combinatorial aspect
connected with the so--called noncrossing partitions \cite{Speicher1997}.

Many questions concerning large random matrices can be easily
reformulated and answered in the framework of operator--valued free
probability; this is also the case of the triangular operator $T$. One can show
that the original definition \eqref{eq:oryginalnadefinicja} is equivalent to
our favourite definition of $T$ as a generalized circular operator which will
be presented in Section \ref{sec:main} \cite{Shlyakhtenko1996,Shlyakhtenko1997}. However,
since in this article we do not need this equivalence, we skip the proof.

In this section we recall briefly some combinatorial aspects the
operator--valued free probability. More information can be found in
papers \cite{Speicher1997,Speicher1998}.

\subsection{Operator--valued probability space}
A triple $(\D\subseteq\A,\E)$ is called an operator--valued probability space
if $\A$ is a unital $\star$--algebra,
$\D\subseteq\A$ is a unital $\star$--subalgebra, and $\E:\A\rightarrow\D$ is a
conditional expectation, i.e.\ $\E$ is linear, satisfies $\E(1)=1$ and
$\E(b_1 x b_2)=b_1 \E(x) b_2$ for every $b_1,b_2\in\D$ and $x\in\A$.


\subsection {Noncrossing pair partitions}
If $X$ is a finite, ordered set, we denote by $\NC_2(X)$ the set of all
noncrossing pair partitions of $X$
\cite{Speicher1998, Kreweras}. A noncrossing pair partition
$\pi=\{Y_1,\dots,Y_n\}$ of $X$ is a decomposition of $X$ into disjoint
two--element sets: $$X=Y_1\cup\cdots\cup Y_n, \qquad Y_i\cap Y_j=
\emptyset \quad \mbox{ if }i\neq j,$$
which has the additional property that for $Y_i=\{a,c\}$ and $Y_j=\{b,d\}$ it
cannot happen that $a<b<c<d$.

We say that the sets $Y_1,\dots,Y_n$ are the lines of the pair partition
$\pi=\{Y_1,\dots,Y_n\}$.
It is useful to describe pair partitions graphically by connecting elements of
the same line by an arc.

\begin{example} There are exactly two noncrossing pair
partitions of the set $\{1,2,3,4\}$ and these are $\big\{ \{1,4\},\{2,3\}
\big\}$ represented by
\begin{pspicture}[](0,-16pt)(0,0)
$\rnode{p1}{1}\ \rnode{p2}{2}\ \rnode{p3}{3}\
\rnode{p4}{4}
\nccurve[angleA=-90,angleB=-90,ncurv=1,nodesep=3pt]{p1}{p4}
\nccurve[angleA=-90,angleB=-90,ncurv=1,nodesep=3pt]{p2}{p3}$ 
\end{pspicture} and
%
%
$\big\{ \{1,2\},\{3,4\} \big\}$ represented by
\begin{pspicture}[](0,-10pt)(0,0) $\rnode{q1}{1}\ \rnode{q2}{2}\ \rnode{q3}{3}\
\rnode{q4}{4}
\nccurve[angleA=-90,angleB=-90,ncurv=1,nodesep=3pt]{q1}{q2}
\nccurve[angleA=-90,angleB=-90,ncurv=1,nodesep=3pt]{q3}{q4}
$
\end{pspicture}.
On the other hand the pair partition $\big\{\{1,3\},\{2,4\} \big\}$ is
crossing as it can be seen on its graphical representation
\begin{pspicture}(0,0)(0,0)
$\rnode{q1}{1}\ \rnode{q2}{2}\ \rnode{q3}{3}\
\rnode{q4}{4}
\nccurve[angleA=-90,angleB=-90,ncurv=1,nodesep=3pt]{q1}{q3}
\nccurve[angleA=-90,angleB=-90,ncurv=1,nodesep=3pt]{q2}{q4}
$
\end{pspicture}.
\end{example}

\subsection{Nested evaluation}
For the purpose of this section we shall forget that $\A$ is an algebra and
that $\D\subseteq\A$. We will assume only that $\D$ is an algebra and that $\A$
is a $\D$--bimodule.

Let there be given a bilinear map $\kk:\A\times\A\rightarrow\D$. We will denote
$$a_1 \bullet_\kk a_2=\kk(a_1,a_2)$$
and assume that $\kk$ is such that
\begin{align}
(b a_1) \bullet_\kk a_2 &=b (a_1 \bullet_\kk a_2), \label{eq:bimodul1} \\
(a_1 b) \bullet_\kk a_2 &=a_1 \bullet_\kk (b a_2), \\
a_1 \bullet_\kk (a_2 b) &=(a_1 \bullet_\kk a_2) b \label{eq:bimodul3}
\end{align}
for any $b\in\D$ and $a_1,a_2\in\A$.


Noncrossing pair partitions of the set $\{1,2,\dots,2n\}$
(or, equivalently, of the set $\{a_1,\dots,a_{2n}\}$) can be identified with
ways of writing brackets in the product $a_1 a_2 \dots a_{2n}$ in such a way
that inside each pair of brackets there are exactly two factors $a_i$ and $a_j$
(and possibly some other nested brackets). To be more explicit: each line
$Y=\{a_i,a_j\}$ of a pair partition corresponds to a certain pair of brackets:
the opening bracket and its closing counterpart. The two factors contained in
between this pair of brackets are exactly $a_i$ and $a_j$.

\begin{example}
$$\begin{pspicture}[](0,-15pt)(0,0) \rnode{q1}{$a_1$}\ \rnode{q2}{$a_2$}\
\rnode{q3}{$a_3$}\ \rnode{q4}{$a_4$}\ \rnode{q5}{$a_5$}\ \rnode{q6}{$a_6$}
\nccurve[angleA=-90,angleB=-90,ncurv=1,nodesep=3pt]{q1}{q2}
\nccurve[angleA=-90,angleB=-90,ncurv=1,nodesep=3pt]{q3}{q4}
\nccurve[angleA=-90,angleB=-90,ncurv=1,nodesep=3pt]{q5}{q6}
\end{pspicture} \text{ corresponds to } (a_1 a_2) (a_3 a_4) (a_5 a_6),$$
$$\begin{pspicture}[](0,-30pt)(0,0)
\rnode{q1}{$a_1$}\ \rnode{q2}{$a_2$}\ \rnode{q3}{$a_3$}\
\rnode{q4}{$a_4$}\ \rnode{q5}{$a_5$}\ \rnode{q6}{$a_6$}
\nccurve[angleA=-90,angleB=-90,ncurv=1,nodesep=3pt]{q1}{q2}
\nccurve[angleA=-90,angleB=-90,ncurv=1,nodesep=3pt]{q3}{q6}
\nccurve[angleA=-90,angleB=-90,ncurv=1,nodesep=3pt]{q4}{q5}
\end{pspicture} \text{ corresponds to } (a_1 a_2) \big( a_3 (a_4 a_5) a_6
\big), $$
$$\begin{pspicture}[](0,-35pt)(0,0)
\rnode{q1}{$a_1$}\ \rnode{q2}{$a_2$}\ \rnode{q3}{$a_3$}\
\rnode{q4}{$a_4$}\ \rnode{q5}{$a_5$}\ \rnode{q6}{$a_6$}
\nccurve[angleA=-90,angleB=-90,ncurv=1,nodesep=3pt]{q1}{q6}
\nccurve[angleA=-90,angleB=-90,ncurv=1,nodesep=3pt]{q2}{q3}
\nccurve[angleA=-90,angleB=-90,ncurv=1,nodesep=3pt]{q4}{q5}
\end{pspicture} \text{ corresponds to } \big(a_1 (a_2 a_3) (a_4 a_5) a_6
\big),$$
$$\begin{pspicture}[](0,-35pt)(0,0)
\rnode{q1}{$a_1$}\ \rnode{q2}{$a_2$}\ \rnode{q3}{$a_3$}\
\rnode{q4}{$a_4$}\ \rnode{q5}{$a_5$}\ \rnode{q6}{$a_6$}
\nccurve[angleA=-90,angleB=-90,ncurv=1,nodesep=3pt]{q1}{q6}
\nccurve[angleA=-90,angleB=-90,ncurv=1,nodesep=3pt]{q2}{q5}
\nccurve[angleA=-90,angleB=-90,ncurv=1,nodesep=3pt]{q3}{q4}
\end{pspicture} \text{ corresponds to } \Big(a_1 \big(a_2 (a_3 a_4) a_5\big)
a_6 \Big).$$
\end{example}

In order to evaluate such a product with brackets we will use the following
rule: in order to multiply two elements of $\A$ we use the product $\bullet_\kk$;
in order to multiply two elements of $\D$ or an element of $\D$ and an element
of $\A$ we use the standard multiplication. In this way for any noncrossing
pair partition $\pi\in\NC_2(\{1,\dots,2n\})$ we have defined a multilinear map
$\kk_\pi:\A^{2n}\rightarrow\D$.

\begin{example}
\begin{multline*} \begin{pspicture}[](0,0)(0,0)
$\kk(\rnode{q1}{a_1},\rnode{q2}{a_2},
\rnode{q3}{a_3},\rnode{q4}{a_4},\rnode{q5}{a_5},\rnode{q6}{a_6})=$
\nccurve[angleA=-90,angleB=-90,ncurv=1,nodesep=3pt]{q1}{q2}
\nccurve[angleA=-90,angleB=-90,ncurv=1,nodesep=3pt]{q3}{q4}
\nccurve[angleA=-90,angleB=-90,ncurv=1,nodesep=3pt]{q5}{q6}
\end{pspicture}
(a_1 \bullet_\kk a_2) (a_3 \bullet_\kk a_4) (a_5 \bullet_\kk a_6)= \\
\kk(a_1,a_2) \kk(a_3,a_4) \kk(a_5,a_6),
\end{multline*}
\begin{multline*}
\begin{pspicture}[](0,0)(0,0)
$\kk(\rnode{q1}{a_1},\rnode{q2}{a_2},\rnode{q3}{a_3},
\rnode{q4}{a_4},\rnode{q5}{a_5},\rnode{q6}{a_6})=$
\nccurve[angleA=-90,angleB=-90,ncurv=1,nodesep=3pt]{q1}{q2}
\nccurve[angleA=-90,angleB=-90,ncurv=1,nodesep=3pt]{q3}{q6}
\nccurve[angleA=-90,angleB=-90,ncurv=1,nodesep=3pt]{q4}{q5}
\end{pspicture}
(a_1 \bullet_\kk a_2) \big(a_3 \bullet_\kk (a_4 \bullet_\kk a_5) a_6\big)=\\
\kk(a_1,a_2) \kk\big( a_3, \kk(a_4,a_5) a_6 \big),
\end{multline*}
\begin{multline*}
\begin{pspicture}[](0,0)(0,0)
$\kk(\rnode{q1}{a_1},\rnode{q2}{a_2},\rnode{q3}{a_3},
\rnode{q4}{a_4},\rnode{q5}{a_5},\rnode{q6}{a_6})=$
\nccurve[angleA=-90,angleB=-90,ncurv=1,nodesep=3pt]{q1}{q6}
\nccurve[angleA=-90,angleB=-90,ncurv=1,nodesep=3pt]{q2}{q3}
\nccurve[angleA=-90,angleB=-90,ncurv=1,nodesep=3pt]{q4}{q5}
\end{pspicture}
a_1 \bullet_\kk (a_2 \bullet_\kk a_3) (a_4 \bullet_\kk a_5) a_6 =\\
\kk\big(a_1, \kk(a_2,a_3)\kk(a_4,a_5) a_6 \big),
\end{multline*}
\begin{multline*}
\begin{pspicture}[](0,0)(0,0)
$\kk(\rnode{q1}{a_1},\rnode{q2}{a_2},\rnode{q3}{a_3},
\rnode{q4}{a_4},\rnode{q5}{a_5},\rnode{q6}{a_6})=$
\nccurve[angleA=-90,angleB=-90,ncurv=1,nodesep=3pt]{q1}{q6}
\nccurve[angleA=-90,angleB=-90,ncurv=1,nodesep=3pt]{q2}{q5}
\nccurve[angleA=-90,angleB=-90,ncurv=1,nodesep=3pt]{q3}{q4}
\end{pspicture}
a_1 \bullet_\kk \big( a_2 \bullet_\kk (a_3 \bullet_\kk a_4) a_5 \big) a_6 = \\
\kk\Big(a_1, \kk\big(a_2, \kk(a_3,a_4) a_5\big) a_6 \Big).
\end{multline*}
\end{example}
\vspace{3ex}

\subsection{Generalized circular elements}
Let an operator--valued probability space $(\D\subseteq\A,\E)$ by given.
We say that $T\in\A$ is a generalized circular element if there exists a
bilinear function $\kk:\A\times\A\rightarrow\D$ (called variance of $T$) which
fulfills \eqref{eq:bimodul1}--\eqref{eq:bimodul3} and such that
\begin{equation}
\E(b_1 T^{s_1} b_2 T^{s_2} \cdots
b_{2n} T^{s_{2n}})=\sum_{\pi\in\NC_2(\{1,\dots,2n\})}
\kk_\pi(b_1 T^{s_1},b_2 T^{s_2},\dots,b_{2n} T^{s_{2n}}),
\label{eq:definicjacyrkularnego}
\end{equation}
\begin{equation}
\E(b_1 T^{s_1} b_2 T^{s_2} \cdots
b_{2n+1} T^{s_{2n+1}})=0
\label{eq:cyrkularnesiezeruja}
\end{equation}
for every $b_1,\dots,b_{2n+1}\in\D$ and $s_1,\dots,s_{2n+1}\in\{1,\star\}$.

\section{The main result}
\label{sec:main}
\subsection{Triangular operator $T$}
Let $\D=\C[x]$ be the $\star$--algebra of polynomials of one variable with
multiplication defined to be the usual multiplication of polynomials and
let $(\D\subset\A,\E)$ be an operator--valued probability space.


The Dykema--Haagerup triangular operator $T\in\A$ is defined to be the generalized
circular element with the variance $\kk$ given by
\begin{equation}
\left\{
\begin{aligned}{}
[\kk(T, b T\gwia)](x)& =\int_x^1 b(t) \dd t, \\ 
[\kk(T\gwia, b T)](x)& =\int_0^x b(t) \dd t, \\ 
[\kk(T, b T)](x)& =0, \\ 
[\kk(T\gwia, b T\gwia)](x)& =0 
\end{aligned}
\right.
\label{eq:deft}
\end{equation}
for any $b\in\D$.

From the following on we shall assume that the algebra $\A$ is generated by the 
algebra $\D$ and operators $T$ and $T\gwia$.

\subsection{Automorphism $\alpha$}
Let us consider a linear automorphism of the algebra $\A$ defined on generators by
$$\alpha(T)=T\gwia, \qquad \alpha(T\gwia)=T,$$
$$[\alpha(b)](x)=b(1-x) \qquad\text{for any } b\in\D.$$
\begin{proposition}
\label{prop:automorfizm}
For any $a\in\A$ we have 
$$\E[\alpha(a)]=\alpha[\E(a)].$$
\end{proposition}
\begin{proof}
It is easy to see that the defining relations
\eqref{eq:deft} give
$$\kk\big( \alpha(T^{s_1}),\alpha(b)
\alpha(T^{s_2})\big)=\alpha\big( \kk(T^{s_1},b T^{s_2}) \big)$$
for any $s_1,s_2\in\{1,\star\}$ and $b\in\D$ and, by induction, that
$$\kk_{\pi}\big( \alpha(b_1 T^{s_1} ),\dots,
\alpha(b_n T^{s_n})\big)=\alpha\big( \kk(b_1 T^{s_1},\dots,b_n T^{s_n}) \big)$$
for any $s_1,\dots,s_n\in\{1,\star\}$, $b_1,\dots,b_n\in\D$ and $\pi\in\NC_2(\{1,\dots,n\})$. 
Finally, we use the definition 
\eqref{eq:definicjacyrkularnego}.
\end{proof}

\subsection{Scalar--valued distribution of $T$}
We define a state $\phi:\A\rightarrow\C$ as follows: for $b\in\D$ we put
$$\phi(b)=\int_0^1 b(x) \dd x,$$
while for a general $a\in\A$ we put
$$\phi(a)=\phi\big( \E(a) \big).$$

\begin{remark} Alternative defintions of the state $\phi$ on the $\star$--algebra 
generated by $T$ 
can be found in articles by Dykema and Haagerup \cite{DykemaHaagerup2001} and 
by Dykema and Yan \cite{DykemaYan}, as well as the proof that $\phi$ is indeed a tracial state. 
\end{remark}

\subsection{Proof of Dykema--Haagerup conjecture}
The following theorem will be our main tool in the proof of Dykema--Haagerup
conjecture. However, since its proof is a bit technical, we postpone it to
Sect. \ref{sec:techniczne}. 
We will use the convention that $T^1=T$ while $T^{-1}=T\gwia$.

\begin{theorem}
\label{theo:formulanapochodna}
Let $s_1,\dots,s_{2m}\in\{1,\star\}$. Then
$\E(T^{s_1} \cdots T^{s_{2m}})$ is a polynomial of degree $m$ and hence
$\frac{d^m}{dx^m} \E(T^{s_1} \cdots T^{s_{2m}})$ can be identified with a 
real number.

Secondly, for any $n\in\N$ we have
\begin{multline}
\label{eq:formulanapochodna}
\frac{d^n}{dx^n} \E(T^{s_1} \cdots T^{s_{2m}})
=\!\!\!\sum_{\substack{0=j_0<j_1<\cdots\\ \cdots<j_{2n}<j_{2n+1}=2m+1}} \\
\left( \frac{d^n}{dx^n} \E(T^{s_{j_1}} T^{s_{j_2}} \cdots T^{s_{j_{2n}}})
\right) \prod_{0\leq r\leq 2m} \E(T^{s_{j_r+1}} T^{s_{j_r+2}} \cdots
T^{s_{j_{r+1}-1}}).
\end{multline}
In the above sum the nonvanishing terms are obtained only for sequences 
$(j_k)$ such that $s_{j_r+1}+s_{j_r+2}+\cdots+s_{j_{r+1}-1}=0$ for every $0\leq r\leq 2m$.

\end{theorem}
\begin{example}
Due to \eqref{eq:cyrkularnesiezeruja} the only nonvanishing terms for $m=3$ and $n=1$ are
\begin{multline*}
\frac{d}{dx} \E (T^{s_1} \cdots T^{s_6})=
\left( \frac{d}{dx} \E (T^{s_1} T^{s_2}) \right)
 \E(T^{s_3} \cdots T^{s_6})+ \\
\left( \frac{d}{dx} \E (T^{s_1} T^{s_4}) \right) \E(T^{s_2} T^{s_3})
\E(T^{s_5} T^{s_6})+
\left( \frac{d}{dx} \E(T^{s_1} T^{s_6}) \right)
\E(T^{s_2} \cdots T^{s_5})+ \\
\left( \frac{d}{dx} \E(T^{s_3} T^{s_4}) \right) \E(T^{s_1} T^{s_2})
\E(T^{s_5} T^{s_6})+
\left( \frac{d}{dx} \E(T^{s_3} T^{s_6}) \right) \E(T^{s_1} T^{s_2})
 \E(T^{s_4} T^{s_5})+\\
\left( \frac{d}{dx} \E(T^{s_5} T^{s_6}) \right) \E(T^{s_1} \cdots T^{s_4}).
\end{multline*}
\end{example}

\begin{theorem}
\label{theo:niezwyklejedynki}
Let $s_1,\dots,s_{2m}\in\{1,-1\}$ be such that
$s_1,s_1+s_2,\dots,s_1+\cdots+s_{2m-1}\leq 0$ and $s_1+\cdots+s_{2m}=0$.
Then
$$\frac{d^m}{dx^m} \E(T^{s_1} \cdots T^{s_{2m}})=1.$$

Let $s_1,\dots,s_{2m}\in\{1,-1\}$ be such that
$s_1,s_1+s_2,\dots,s_1+\cdots+s_{2m-1}\geq 0$ and $s_1+\cdots+s_{2m}=0$.
Then
$$\frac{d^m}{dx^m} \E(T^{s_1} \cdots T^{s_{2m}})=(-1)^m.$$
\end{theorem}
\begin{proof}
We shall prove the first part of the theorem by induction (the proof of the
second part is analogous and we skip it).

Let us compute $\frac{d^{m-1}}{dx^{m-1}} \E(T^{s_1} \cdots T^{s_{2m}})$
from Theorem \ref{theo:formulanapochodna}.
It is easy to observe that it yields that
\begin{multline*}
\frac{d^{m-1}}{dx^{m-1}} \E(T^{s_1} \cdots T^{s_{2m}})= \\
\sum_{1\leq k\leq 2m-1}
\left( \frac{d^{m-1}}{dx^{m-1}} \E(T^{s_1} T^{s_2}
\cdots T^{s_{k-1}} T^{s_{k+2}} T^{s_{k+3}} \cdots T^{s_{2m}}) \right)
\E(T^{s_k} T^{s_{k+1}}).
\end{multline*}
For nonzero summands the inductive hypothesis can be applied and hence
$$\frac{d^{m-1}}{dx^{m-1}} \E(T^{s_1} \cdots T^{s_{2m}})= \\
\sum_{1\leq k\leq 2m-1}
\E(T^{s_k} T^{s_{k+1}}).$$

From \eqref{eq:deft} follows that if $s_k=s_{k+1}$ then $\E(T^{s_k}
T^{s_{k+1}})=0$. Since the sequence $s_1,\dots,s_{2m}$ must begin with $-1$ and end with
$1$, hence for some $l\geq 0$ there are exactly $l$ values of the index $1\leq k<2m$
such that $(s_{k},s_{k+1})=(1,-1)$ and exactly $l+1$ values of the index $1\leq k<2m$
such that $(s_{k},s_{k+1})=(-1,1)$. Therefore
$$\frac{d^{m}}{dx^{m}} \E(T^{s_1} \cdots T^{s_{2m}})=
l \frac{d}{dx} \E(T^{1} T^{-1})+(l+1) \frac{d}{dx} \E(T^{-1} T^{1})=1.$$

\end{proof}
\begin{remark}
We leave to the reader the proof of the following: suppose that for $m\in\N$ we
have $s_1,\dots,s_{2m},s_1',\dots,s_{2m}'\in\{1,-1\}$ such that
$s_1+\cdots+s_{2m}=s_1'+\cdots+s_{2m}'=0$ and that for every $1\leq k\leq 2m$
the sums $s_1+\cdots+s_k$ and $s'_1+\cdots+s'_k$ have the same sign (to be precise:
$(s_1+\cdots+s_k)(s'_1+\cdots+s'_k)\geq 0$). Then
$$\frac{d^m}{dx^m} \E(T^{s_1} \cdots T^{s_{2m}})=
\frac{d^m}{dx^m} \E(T^{s_1'} \cdots T^{s_{2m}'}).$$
 \end{remark}

\begin{proposition}
\label{prop:nihilizm}
Let $s_1,\dots,s_{2m}\in\{1,\star\}$. If $s_1=1$ or $s_{2m}=\star$ then
$$\E(T^{s_1} \cdots T^{s_{2m}})\big|_{x=1}=0.$$
\end{proposition}
\begin{proof}
From defining relations \eqref{eq:deft} we obtain that for any $b\in\D$
and $s\in\{1,\star\}$ we have $\E(T b T^{s})|_{x=1}=\E(T^{s} 
T\gwia)|_{x=1}=0$. Now observe that in the evaluation of a nested product
$\kk_{\pi}(T^{s_1},\dots,T^{s_{2m}})$ one of the above expressions must appear
(if observations are not successful we refer to
Eq.\ \eqref{eq:rekursjanacyrkularne}).
\end{proof}

\begin{theorem}
For every $k,n\in\N$, $x\in\R$ and $0\leq m\leq k-1$
\label{theo:indukcja}
\begin{align}
\label{eq:poczatkowe}
\frac{d^m}{dx^m} \E\big[ \big(T^k (T\gwia)^k \big)^n \big]
\bigg|_{x=1}& =0, \\ 
\label{eq:rekurencja}
\frac{d^k}{dx^k} \E\big[ \big(T^k (T\gwia)^k \big)^n \big](x) & =(-1)^k
\E\big[ \big( T^k (T\gwia)^k\big)^{n-1} \big] (x-1).
\end{align}
\end{theorem}
\begin{remark}
Equations \eqref{eq:poczatkowe} and \eqref{eq:rekurencja} allow us to
find $\E\big[ \big(T^k (T\gwia)^k \big)^n \big]$ uniquely if
$\E\big[ \big(T^k (T\gwia)^k \big)^{n-1} \big]$ is known.
In particular, for $k=1$ equation \eqref{eq:rekurencja} coincides with the 
defining relation of Abel polynomials (up to normalisation), hence
$$ \E\big[ (T \ T\gwia)^n \big](1-x) = \E\big[ (T\gwia T)^n \big](x)=
\\ \frac{1}{n!} A_n(x) = \frac{x (x-n)^{n-1}}{n!}. $$ 
\end{remark}
\begin{proof}
We use Theorem \ref{theo:formulanapochodna} in order to compute
$\frac{d^m}{dx^m}\E\big[ \big(T^k (T\gwia)^k \big)^n \big] \big|_{x=1}$.
Proposition \ref{prop:nihilizm} implies that nonzero terms will be obtained
only for sequences
\begin{multline*}
j_1=1, j_2=2, \dots, j_k=k< j_{k+1} < j_{k+2} < \cdots <
j_{2m-k}< \\
<j_{2m-k+1}=2kn-k+1,j_{2m-k+2}=2kn-k+2,\dots,
j_{2m}=2kn.
\end{multline*}
It has twofold implications. Firstly, since the above condition cannot be
fulfilled for $m<k$, the first part of the theorem follows.

Secondly, for every $m\geq 0$
$$\frac{d^{m+k}}{dx^{m+k}}\E\big[ \big(T^k (T\gwia)^k \big)^n \big]
\bigg|_{x=1}=
(-1)^{k+m} \frac{d^m}{dx^m}\E\big[ \big((T\gwia)^k T^k\big)^{n-1} \big]
\bigg|_{x=1}$$
by applying Eq. \eqref{eq:formulanapochodna} and 
Theorem \ref{theo:niezwyklejedynki} to both sides of the equation (we leave it as
an exercise to the reader to verify that assumptions of Theorem \ref{theo:niezwyklejedynki}
are fulfilled for nonvanishing summands). 
It follows that
$$\frac{d^{k}}{dx^{k}}\E\big[ \big(T^k (T\gwia)^k \big)^n \big](1+x)=
(-1)^k \E\big[ \big((T\gwia)^k T^k\big)^{n-1} \big](1-x)$$
since both sides of the above equation are polynomials, all derivatives of
which coincide in $x=0$.

Now observe that
\begin{multline*}
(-1)^k \E\big[ \big((T\gwia)^k T^k\big)^{n-1} \big](1-x) =
(-1)^k \E\Big[\alpha\Big( \big((T\gwia)^k T^k\big)^{n-1} \Big)\Bigg](x)=\\
(-1)^k \E\big[ \big(T^k (T\gwia)^k\big)^{n-1} \big](x)
\end{multline*}
which finishes the proof.
\end{proof}
\begin{remark}
Let $k_1\geq k_2 \geq \cdots \geq k_n\geq 1$. By
$T^{k_1} (T\gwia)^{k_2} T^{k_3} \cdots$ we denote the alternating product of powers of
$T$ and $T\gwia$ the last factor of which is equal to $T^{k_n}$ (for $n$ odd) or equal to
$(T\gwia)^{k_n}$ (for $n$ even). We leave to the reader the proof of the following:
if $0\leq m\leq k_1-1$ then
$$\frac{d^m}{dx^m} \E\big[\big( T^{k_1} (T\gwia)^{k_2} T^{k_3} \cdots \big) \big( \cdots
(T\gwia)^{k_3} T^{k_2} (T\gwia)^{k_1} \big)
\big]\bigg|_{x=1} =0, $$
\begin{multline*}
\frac{d^{k_1}}{dx^{k_1}} 
\E\big[ \big(T^{k_1} (T\gwia)^{k_2} T^{k_3} \cdots\big) 
\big(\cdots (T\gwia)^{k_3} T^{k_2} (T\gwia)^{k_1}
\big) \big](x)= \\
(-1)^{k_1} \E\big[\big(T^{k_2} (T\gwia)^{k_3} T^{k_4} \cdots\big)
\big(\cdots (T\gwia)^{k_4} T^{k_3} (T\gwia)^{k_2} \big)\big](x-1).
\end{multline*}

\end{remark}

\begin{corollary}
\label{coro:objetosc}
Let $k,n\in\N$. For every $x<1$ we have
\begin{multline}
\label{eq:objetosc1}
\E\big[ \big(T^k (T\gwia)^k\big)^n \big](x)=\vol \Big\{
(x_1,\dots,x_{nk})\in\R^{nk}: \\ x<x_1<x_2<\cdots<x_{kn} \text{ and }
x_{k}<1 \text{ and } x_{2k}<2,\cdots,\text{ and } x_{nk}<n \Big\}.
\end{multline}
Furthermore
$$\phi\big[ \big(T^k (T\gwia)^k\big)^n \big]=\vol V_{k,n},$$
where
\begin{multline*} V_{k,n}=\Big\{
(x_0,x_1,\dots,x_{nk})\in\R^{nk+1}:  0<x_0<x_1<\cdots<x_{kn} \text{ and }\\
x_{k}<1 \text{ and } x_{2k}<2,\cdots,\text{ and } x_{nk}<n \Big\}.
\end{multline*}
\end{corollary}
\begin{proof}
From Theorem \ref{theo:indukcja} it follows that
\begin{multline*}
\E\big[ \big(T^k (T\gwia)^k\big)^n \big](x)=    \\
\int_x^1 \int_{x_1}^1 \cdots \int_{x_{k-1}}^1
\E\big[ \big( T^k (T\gwia)^k\big)^{n-1} \big](x_k-1) \dd x_1 \dd x_2
\cdots \mathrm{d} x_k
\end{multline*}
and since the volume in (\ref{eq:objetosc1}) can be written as an iterated
integral, the first part of the Lemma follows by induction.

The second part of the Corollary is a direct consequence of the first one.
\end{proof}

\begin{theorem}[Conjecture of Dykema and Haagerup]
For any $k,n\in\N$ we have
$$\phi\big[ \big(T^k (T\gwia)^k\big)^n \big]=\frac{n^{nk}}{(nk+1)!}.$$
\end{theorem}
\begin{proof}
We denote
\begin{multline*}
S_{k,n}=\big\{ \x=(x_0,\dots,x_{nk})\in\R^{nk+1}:  \\
0<x_0<\cdots<x_{nk}<n, \qquad x_0,\dots,x_{nk}\not\in\Z \big\}
\end{multline*}
and for $l\in\{0,1,\dots,n-1\}$ and $\x\in S_{k,n}$ we define
$s_l(\x)+\frac{nk+1}{n}$ to be the number of elements of the tuple $\x$ which
are also elements of the interval $(l,l+1)$:
$$s_l(\x)=\overline{\overline{\{x_0,\dots,x_{nk}\}\cap [l,l+1]
}}-\frac{nk+1}{n}.$$
So defined numbers fulfill
$$s_0(\x)+\cdots+s_{n-1}(\x)=0.$$
For $\x\in S_{k,n}$ it is easy to see that $\x\in V_{k,n}$ if and only
if
$$s_0(\x),s_0(\x)+s_1(\x),\dots,s_0(\x)+\cdots+s_{n-2}(\x)>0.$$

For any integer number $m$ and $\x\in S_{k,n}$ we define $m+\x\in S_{k,n}$ to
be the increasing sequence such that the set of values of $m+\x$ is equal to
$$\{m+x_0 \bmod n, m+x_1 \bmod n,\dots,m+x_{nk} \bmod n\}.$$
We recall that if $y\in\R$ then $y \bmod n$ is defined to be the $y'$ such
that $y'-y$ is a multiple of $n$ and $0\leq y'<n$.

Note that $s_l(m+\x)=s_{(l-m) \bmod n}(\x)$ and hence the sequences
$$\big(s_0(m+\x),\dots,s_{n-1}(m+\x)\big)_{m=0,1,2,\dots,n-1}$$
are cyclic rotations of each other.

From Raney lemma \cite{Raney1960} it
follows that for every $\x\in S_{k,n}$ there exists exactly one
$m\in\{0,1,\dots,n-1\}$ such that $m+\x\in V_{k,n}$; in order to be self--contained
we indicate the proof. Indeed, it is the $m$ for
which the sum $s_1(\x)+\cdots+s_{m}(\x)$ takes its minimal value.
Such index $m$ is unique, since if
$m_1\neq m_2$ then
$$\big(s_1(\x)+\cdots+s_{m_1}(\x)\big)-\big(s_1(\x)+\cdots+s_{m_2}(\x)\big)$$
cannot be integer and hence
$$s_1(\x)+\cdots+s_{m_1}(\x)\neq s_1(\x)+\cdots+s_{m_2}(\x).$$

We have proved that $S_{k,n}$ is a disjoint sum of the sets 
$(m+V_{k,n})_{m=0,1,\dots,n-1}$. Since all these sets have the same volume, we have
$$n \vol V_{k,n}=\frac{n^{nk+1}}{(nk+1)!}$$ which finishes the proof after application of
Corollary \ref{coro:objetosc}.
\end{proof}

\begin{remark}
In this article we define $\D$ to be an algebra of polynomials on the whole real line.
From the viewpoint of the theory of large random matrices developed by
Shlyakhtenko \cite{Shlyakhtenko1996,Shlyakhtenko1997} it does not have much sense, since
only the values of functions on the interval $[0,1]$ have a nice interpretation. Nevertheless,
our proof of Dykema--Haagerup conjecture uses extensively such senseless objects. It
is a small mystery that we still do not understand. 
\end{remark}

\begin{remark}
By refining the above proof one can show the following: for every sequence 
$l_1\geq l_2 \geq \cdots \geq l_n$ it is possible to find a close formula for a
family of moments
$$\phi\big[ \big(T^{k+l_1} (T\gwia)^{k+l_2} T^{k+l_3} \cdots\big) 
\big(\cdots (T\gwia)^{k+l_3} T^{k+l_2} (T\gwia)^{k+l_1}
\big) \big].$$

For example,
$$\phi\big( T^{n+1} (T\gwia)^n T^n (T\gwia)^{n+1} \big)=\frac{1}{2} \left( 
\frac{2^{2n+2}}{(2n+2)!} - \frac{1}{[(n+1)!]^2} \right), $$
$$\phi\big( T^{n+1} (T\gwia)^n T^n (T\gwia)^n T^n (T\gwia)^{n+1} \big)=
\frac{3^{3n+1}}{(3n+2)!}- \frac{1}{n! [(n+1)!]^2 }, $$
\begin{multline*}
\phi\big( T^{n+2} (T\gwia)^n T^n (T\gwia)^n T^n (T\gwia)^{n+2} \big)= \\
\frac{3^{3n+2}}{(3n+3)!}- 
\frac{2^{2n}}{(n+2)! (2n+1)!} - \frac{2^{2n+3}}{3 (n+1)! (2n+2)!}
+\frac{1}{3 [(n+1)!]^3}.
\end{multline*}
We leave to the reader to apply methods from Sect.\ \ref{sec:kombinatoryczna}
to the left hand sides of the above equations
and obtain even more multlinomial identities 
similar to \eqref{eq:tozsamosc2}--\eqref{eq:tozsamosc4}.
\end{remark}

\section{Formulas for $\phi\big[ \big(T^k (T\gwia)^k\big)^n \big]$ involving multinomial 
coefficients}
\label{sec:kombinatoryczna}

For any fixed $n$ it is possible to write a (possibly very complicated) formula for
$\phi\big[ \big(T^k (T\gwia)^k\big)^n \big]$ which involves multinomial coefficients. 
The first step is to enumerate all noncrossing pair partitions of the tuple 
\begin{equation}
\underbrace{\underbrace{T,\dots,T}_{k \text{ times}},\underbrace{T\gwia,\dots,T\gwia}_{k 
\text{ times}},
\dots,\underbrace{T,\dots,T}_{k \text{ times}},
\underbrace{T\gwia,\dots,T\gwia}_{k \text{ times}}}_{2n \text{ times}}
\label{eq:mazursniady}
\end{equation}
with a property that each line connects either $T$ with $T\gwia$ or $T\gwia$ with $T$ 
(see the example below).

Secondly, for $\alpha,\beta\geq 0$ we define a function $f_{\alpha,\beta}\in\D$ by
$$f_{\alpha,\beta}(x)=\frac{x^\alpha (1-x)^\beta}{\alpha! \beta!}.$$
It is not difficult to show the following lemma.
\begin{lemma}
\label{lem:rachownictwo}
For every $\alpha,\beta,\alpha',\beta',p\geq 0$ we have
$$f_{\alpha,\beta} f_{\alpha',\beta'}=\binom{\alpha+\alpha'}{\alpha}
 \binom{\beta+\beta'}{\beta} f_{\alpha+\alpha',\beta+\beta'},$$
$$\begin{pspicture}[](0,-20pt)(0,20pt) 
$\kk(\overbrace{ \rnode{r1}{T}}^{p \text{ lines}},f_{\alpha,\beta} 
\rnode{s1}{T\gwia})=$
\ncangle[angleA=-90,angleB=-90,ncurv=1,arm=10pt,nodesep=1pt,linearc=10pt,doubleline=true]{r1}{s1}
\end{pspicture}
\sum_{0\leq k\leq \alpha} \binom{k+p}{p} f_{\alpha-k,\beta+k+p},$$
$$\begin{pspicture}[](0,-20pt)(0,0) 
$\kk(\overbrace{ \rnode{r1}{T\gwia}}^{p \text{ lines}},f_{\alpha,\beta} 
\rnode{s1}{T})=$
\ncangle[angleA=-90,angleB=-90,ncurv=1,arm=10pt,nodesep=1pt,linearc=10pt,doubleline=true]{r1}{s1}
\end{pspicture}
\sum_{0\leq k\leq \alpha} \binom{k+p}{p} f_{\alpha+k+p,\beta-k},$$
$$\begin{pspicture}[](0,-20pt)(0,20pt) 
$\phi\big(\kk(\overbrace{ \rnode{t1}{T}}^{p \text{ lines}},f_{\alpha,\beta} 
\rnode{u1}{T\gwia})\big)
=\frac{\binom{p+\alpha}{\alpha}}{(\alpha+\beta+p+1)!},$ 
\ncangle[angleA=-90,angleB=-90,ncurv=1,arm=10pt,nodesep=1pt,linearc=10pt,doubleline=true]{r1}{s1}
\ncangle[angleA=-90,angleB=-90,ncurv=1,arm=10pt,nodesep=1pt,linearc=10pt,doubleline=true]{t1}{u1}
\end{pspicture} $$
$$\begin{pspicture}[](0,-20pt)(0,20pt) 
$\phi\big(\kk(\overbrace{ \rnode{t1}{T\gwia}}^{p \text{ lines}},f_{\alpha,\beta} 
\rnode{u1}{T})\big)
=\frac{\binom{p+\beta}{\beta}}{(\alpha+\beta+p+1)!},$ 
\ncangle[angleA=-90,angleB=-90,ncurv=1,arm=10pt,nodesep=1pt,linearc=10pt,doubleline=true]{r1}{s1}
\ncangle[angleA=-90,angleB=-90,ncurv=1,arm=10pt,nodesep=1pt,linearc=10pt,doubleline=true]{t1}{u1}
\end{pspicture} $$
where in order to save space and keep the notation simple, instead of
$$\begin{pspicture}[](0,-20pt)(0,20pt) 
$\overbrace{ \rnode{a1}{T},\dots,\rnode{a3}{T}}^{p \text{ lines}},$ \hspace{2cm} $
\rnode{x3}{T\gwia},\dots,\rnode{x1}{T\gwia}
\ncangle[angleA=-90,angleB=-90,ncurv=1,arm=15pt,nodesep=1pt,linearc=15pt]{a1}{x1}
\ncangle[angleA=-90,angleB=-90,ncurv=1,arm=10pt,nodesep=1pt,linearc=10pt]{a3}{x3}
\text{ we write }
\overbrace{ \rnode{r1}{T}}^{p \text{ lines}},$ \hspace{2cm}
$\rnode{s1}{T\gwia}.$
\ncangle[angleA=-90,angleB=-90,ncurv=1,arm=10pt,nodesep=1pt,linearc=10pt,doubleline=true]{r1}{s1}
\end{pspicture}$$
\end{lemma}

The third and last step to find the formula for $\phi\big[ \big(T^k (T\gwia)^k\big)^n \big]$ 
is a direct application of Lemma \ref{lem:rachownictwo} to compute the nested product 
$\kk_\pi \big[ \big(T^k (T\gwia)^k\big)^n \big]$.

The last two steps can also be replaced by a more combinatorial procedure of counting 
orderings on certain graphs \cite{DykemaHaagerup2001}.

\begin{example}
We have that for $n=3$ every noncrossing pair partition (with the additional property 
mentioned above) of the tuple \eqref{eq:mazursniady} is in the form 
\begin{equation}
\label{eq:postac0}
\begin{pspicture}[](0,-25pt)(0,25pt)
$
\overbrace{ \rnode{a1}{T}}^{p \text{ lines}},
\overbrace{ \rnode{q1}{T}}^{q \text{ lines}},\rnode{r1}{T\gwia},
\overbrace{\rnode{c1}{T\gwia}}^{p \text{ lines}},
\rnode{d1}{T},
\overbrace{ \rnode{j1}{T}}^{q \text{ lines}},
\rnode{k1}{T\gwia},
\overbrace{\rnode{h1}{T\gwia}}^{p \text{ lines}},
\rnode{i1}{T},
\overbrace{\rnode{e1}{T}}^{q \text{ lines}},
\rnode{f1}{T\gwia},\rnode{b1}{T\gwia},$
\ncangle[angleA=-90,angleB=-90,ncurv=1,arm=20pt,nodesep=1pt,linearc=20pt,doubleline=true]{a1}{b1}
\ncangle[angleA=-90,angleB=-90,ncurv=1,arm=10pt,nodesep=1pt,linearc=10pt,doubleline=true]{r1}{q1}
\ncangle[angleA=-90,angleB=-90,ncurv=1,arm=10pt,nodesep=1pt,linearc=10pt,doubleline=true]{c1}{d1}
\ncangle[angleA=-90,angleB=-90,ncurv=1,arm=10pt,nodesep=1pt,linearc=10pt,doubleline=true]{h1}{i1}
\ncangle[angleA=-90,angleB=-90,ncurv=1,arm=10pt,nodesep=1pt,linearc=10pt,doubleline=true]{j1}{k1}
\ncangle[angleA=-90,angleB=-90,ncurv=1,arm=10pt,nodesep=1pt,linearc=10pt,doubleline=true]{e1}{f1} 
\end{pspicture}
\end{equation}
with $p+q=k$ or in one of the forms
\begin{equation}
\begin{pspicture}[](0,-40pt)(0,25pt)
$
\overbrace{ \rnode{a1}{T}}^{p \text{ lines}},
\overbrace{ \rnode{x1}{T}}^{q+1 \text{ lines}},
\overbrace{ \rnode{q1}{T}}^{r \text{ lines}},
\rnode{r1}{T\gwia},\!\!\!\!
\overbrace{\rnode{c1}{T\gwia}}^{p+q+1 \text{ lines}}\!\!\!\!\,
\rnode{d1}{T},\overbrace{ \rnode{j1}{T}}^{r \text{ lines}},
\rnode{k1}{T\gwia},\rnode{y1}{T\gwia},
\overbrace{\rnode{h1}{T\gwia}}^{p \text{ lines}},
\rnode{i1}{T},\!\!\!\!\overbrace{ \rnode{e1}{T}}^{q+r+1 \text{ lines}}\!\!\!\!,
\rnode{f1}{T\gwia},\rnode{b1}{T\gwia},$
\ncangle[angleA=-90,angleB=-90,ncurv=1,arm=30pt,nodesep=1pt,linearc=30pt,doubleline=true]{a1}{b1}
\ncangle[angleA=-90,angleB=-90,ncurv=1,arm=20pt,nodesep=1pt,linearc=20pt,doubleline=true]{x1}{y1}
\ncangle[angleA=-90,angleB=-90,ncurv=1,arm=10pt,nodesep=1pt,linearc=10pt,doubleline=true]{r1}{q1}
\ncangle[angleA=-90,angleB=-90,ncurv=1,arm=10pt,nodesep=1pt,linearc=10pt,doubleline=true]{c1}{d1}
\ncangle[angleA=-90,angleB=-90,ncurv=1,arm=10pt,nodesep=1pt,linearc=10pt,doubleline=true]{h1}{i1}
 \ncangle[angleA=-90,angleB=-90,ncurv=1,arm=10pt,nodesep=1pt,linearc=10pt,doubleline=true]{j1}{k1}
\ncangle[angleA=-90,angleB=-90,ncurv=1,arm=10pt,nodesep=1pt,linearc=10pt,doubleline=true]{e1}{f1}
\end{pspicture}
\end{equation}
\begin{equation}
\begin{pspicture}[](0,-40pt)(0,25pt)
$
\!\!\!\!\overbrace{ \rnode{a1}{T}}^{p+q+1 \text{ lines}},
\overbrace{ \rnode{q1}{T}}^{r \text{ lines}},
\rnode{r1}{T\gwia},
\overbrace{ \rnode{x1}{T\gwia}}^{q \text{ lines}},
\overbrace{\rnode{c1}{T\gwia}}^{p \text{ lines}},
\rnode{d1}{T},
\!\!\!\!\overbrace{ \rnode{j1}{T}}^{q+r+1 \text{ lines}}\!\!\!\!,
\rnode{k1}{T\gwia},
\overbrace{ \rnode{h1}{T\gwia}
}^{p \text{ lines}},
\rnode{i1}{T},\rnode{y1}{T},
\overbrace{ \rnode{e1}{T}}^{r \text{ lines}},
\rnode{f1}{T\gwia},\rnode{b1}{T\gwia},$
\ncangle[angleA=-90,angleB=-90,ncurv=1,arm=30pt,nodesep=1pt,linearc=30pt,doubleline=true]{a1}{b1}
\ncangle[angleA=-90,angleB=-90,ncurv=1,arm=10pt,nodesep=1pt,linearc=10pt,doubleline=true]{r1}{q1}
\ncangle[angleA=-90,angleB=-90,ncurv=1,arm=20pt,nodesep=1pt,linearc=20pt,doubleline=true]{x1}{y1}
\ncangle[angleA=-90,angleB=-90,ncurv=1,arm=10pt,nodesep=1pt,linearc=10pt,doubleline=true]{c1}{d1}
\ncangle[angleA=-90,angleB=-90,ncurv=1,arm=10pt,nodesep=1pt,linearc=10pt,doubleline=true]{h1}{i1}
\ncangle[angleA=-90,angleB=-90,ncurv=1,arm=10pt,nodesep=1pt,linearc=10pt,doubleline=true]{j1}{k1}
\ncangle[angleA=-90,angleB=-90,ncurv=1,arm=10pt,nodesep=1pt,linearc=10pt,doubleline=true]{e1}{f1}
\end{pspicture}
\end{equation}
\begin{equation}
\label{eq:postac3}
\begin{pspicture}[](0,-40pt)(0,25pt)
$
\overbrace{ \rnode{a1}{T}}^{p \text{ lines}},
\overbrace{ \rnode{q1}{T}}^{q+r+1 \text{ lines}}\!\!\!\!,
\rnode{r1}{T\gwia},
\overbrace{ \rnode{c1}{T\gwia}}^{p \text{ lines}},
\rnode{d1}{T},\overbrace{ \rnode{x1}{T}}^{q \text{ lines}},
\overbrace{ \rnode{j1}{T}}^{r \text{ lines}},
\rnode{k1}{T\gwia},
\!\!\!\!\overbrace{ \rnode{h1}{T\gwia}}^{p+q+1 \text{ lines}}\!\!\!\!,
\rnode{i1}{T},\overbrace{ \rnode{e1}{T}}^{r \text{ lines}},
\rnode{f1}{T\gwia},\rnode{y1}{T\gwia},\rnode{b1}{T\gwia}$
\ncangle[angleA=-90,angleB=-90,ncurv=1,arm=30pt,nodesep=1pt,linearc=30pt,doubleline=true]{a1}{b1}
\ncangle[angleA=-90,angleB=-90,ncurv=1,arm=10pt,nodesep=1pt,linearc=10pt,doubleline=true]{r1}{q1}
\ncangle[angleA=-90,angleB=-90,ncurv=1,arm=10pt,nodesep=1pt,linearc=10pt,doubleline=true]{c1}{d1}
\ncangle[angleA=-90,angleB=-90,ncurv=1,arm=20pt,nodesep=1pt,linearc=20pt,doubleline=true]{x1}{y1}
\ncangle[angleA=-90,angleB=-90,ncurv=1,arm=10pt,nodesep=1pt,linearc=10pt,doubleline=true]{h1}{i1}
\ncangle[angleA=-90,angleB=-90,ncurv=1,arm=10pt,nodesep=1pt,linearc=10pt,doubleline=true]{j1}{k1}
\ncangle[angleA=-90,angleB=-90,ncurv=1,arm=10pt,nodesep=1pt,linearc=10pt,doubleline=true]{e1}{f1}
\end{pspicture}
\end{equation}
with $p+q+r=k-1$.

We leave it as an exercise to the reader to apply Lemma \ref{lem:rachownictwo} and
check that $\phi\big[ \big(T^k (T\gwia)^k\big)^3 \big]$ is indeed equal to the 
right hand side of \eqref{eq:tozsamosc3}. 
\end{example}

\section{Technical results}
\label{sec:techniczne}

The main result of this Section is the proof of Theorem
\ref{theo:formulanapochodna}. 

\begin{proposition}
\label{prop:innadefinicja}
If $T$ is a generalized circular element, $m\geq 0$ and
$s_1,\dots,s_{2m}\in\{1,\star\}$ then
\begin{multline}
\E(T^{s_1} \cdots T^{s_{2m}})=
\sum_{k\geq 1}
\sum_{\substack{1=i_1<i_2<\cdots\\ \cdots<i_{k+1}=2m+1}}
\prod_{1\leq r\leq k}
\label{eq:rekursjanacyrkularne}  \\
\kk\big( T^{s_{i_r}}, \E(T^{s_{i_r+1}} T^{s_{i_r+2}} \cdots
T^{s_{i_{r+1}-3}} T^{s_{i_{r+1}-2}}) T^{s_{i_{r+1}-1}} \big).
\end{multline}
\end{proposition}
\begin{remark}
We can treat Proposition \ref{prop:innadefinicja} as an alternative recursive
definition of the generalized circular elements.
\end{remark}
\begin{example}
If $T$ is a generalized circular element then for any
$s_1,\dots,s_6\in\{1,\star\}$ we have
\begin{multline*}
\E(T^{s_1} T^{s_2} T^{s_3} T^{s_4} T^{s_5} T^{s_6})= \\
\kk\big(T^{s_1}, \E(T^{s_2} T^{s_3} T^{s_4} T^{s_5}) T^{s_6}\big)+
\kk\big(T^{s_1}, T^{s_2}) \kk \big(T^{s_3},
                            \E( T^{s_4} T^{s_5}) T^{s_6}\big)+\\
\kk\big(T^{s_1}, \E(T^{s_2} T^{s_3}) T^{s_4}\big) \kk(T^{s_5},T^{s_6})+
\kk(T^{s_1}, T^{s_2}) \kk(T^{s_3}, T^{s_4}) \kk(T^{s_5},T^{s_6}).
\end{multline*}
\end{example}
\begin{proof}
%
Observe that for $\pi\in\NC_2(\{1,\dots,2m\})$ we can distinguish certain lines
in $\pi$ which will be called outer lines. This name can be easily justified by
the example below, where the outer lines were plotted with a bold line.
To be precise: we define $i_1=1$ and we inductively define $i_{n}$ by the
requirement that the index $i_{n}-1$ is joined by a line in $\pi$ with
the index $i_{n-1}$.
In this way we have defined a tuple $1=i_1<i_2<\cdots<i_{k+1}=2m+1$ and the
outer lines in $\pi$ are exacly $\{i_1,i_2-1\}, \{i_2,i_3-1\}, \dots,
\{i_{k},i_{k+1}-1\}$.

Furthermore, for $1\leq r\leq k$ we define
$\pi_r\in\NC_2(\{i_{r}+1,i_{r}+2,\dots,i_{r+1}-2\})$ by $\pi_r=\big\{
\{a,b\}\in\pi: i_{r}+1\leq a,b \leq i_{r+1}-2 \big\}$. Lines of noncrossing
pair partitions $\pi_1,\dots,\pi_k$ have a nice graphical interpretation as
inner lines of the partition $\pi$ \cite{BozejkoLeinertSpeicher}.

\begin{example}
For a noncrossing pair partition $\pi$ given by
$$
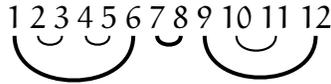
\begin{pspicture}[](0,-20pt)(0,0)
$\rnode{q1}{1}\ \rnode{q2}{2}\ \rnode{q3}{3}\
\rnode{q4}{4}\ \rnode{q5}{5}\ \rnode{q6}{6}\ \rnode{q7}{7}\ \rnode{q8}{8}
\ \rnode{q9}{9}\ \rnode{q10}{10}\ \rnode{q11}{11}\ \rnode{q12}{12}$
\nccurve[linewidth=1.5pt,angleA=-90,angleB=-90,
ncurv=1,nodesep=3pt]{q1}{q6}
\nccurve[angleA=-90,angleB=-90,ncurv=1,nodesep=3pt]{q2}{q3}
\nccurve[angleA=-90,angleB=-90,ncurv=1,nodesep=3pt]{q4}{q5}
\nccurve[linewidth=1.5pt,angleA=-90,angleB=-90,ncurv=1,
 nodesep=3pt]{q7}{q8}
\nccurve[linewidth=1.5pt,angleA=-90,angleB=-90,ncurv=1,
 nodesep=3pt]{q9}{q12}
\nccurve[angleA=-90,angleB=-90,ncurv=1,
 nodesep=3pt]{q10}{q11}
\end{pspicture}
$$
we have $i_1=1$, $i_2=7$, $i_3=9$, $i_4=13$. The outer lines are: $\{1,6\}$,
$\{7,8\}$, $\{9,12\}$ and are drawn with a bold line. We have
$\pi_1=\big\{ \{2,3\},\{4,5\} \big\}$, $\pi_2=\emptyset$,
$\pi_3=\big\{ \{10,11\} \big\}$.
\end{example}

Now we see that the sum over $\pi$ in \eqref{eq:definicjacyrkularnego} can be
replaced by the sum over $k\geq 1$ (i.e.\ the number of outer lines of $\pi$),
over  indices $1=i_1<i_2<\cdots<i_{k+1}=2m+1$ (i.e.\ the positions of the outer
lines of $\pi$) and over $\pi_1,\dots,\pi_k$ such that
$\pi_r\in\NC_2(\{i_{r}+1,i_{r}+2,\dots,i_{r+1}-2\})$ (i.e.\ the inner
lines of $\pi$):

\begin{multline*}
\E(T^{s_1}\dots T^{s_{2m}})=\\ \shoveleft{\sum_{k\geq 1}
\sum_{i_1<\cdots<i_{k+1}}
\begin{pspicture}[](0,0)(0,-15pt) $
\kk\Big(\rnode{p1}{T^{s_{i_1}}}, \sum_{\pi_1}
\kk_{\pi_1}(T^{s_{i_1+1}},\dots,T^{s_{i_2
-2}}) \rnode{p2}{T^{s_{i_2-1}}},
\ncangle[angleA=-90,angleB=-90,ncurv=1,nodesep=3pt,linearc=.30]{p1}{p2}
$ \end{pspicture}} \\
\dots, \begin{pspicture}[](0,0)(0,0) $
\rnode{p3}{T^{s_{i_{k}}}}, \sum_{\pi_k}
\kk_{\pi_k}(T^{s_{i_{k}+1}},\dots,T^{s_{i_{k+1}
-2}})\rnode{p4}{T^{s_{i_{k+1}-1}}}\Big).
\ncangle[angleA=-90,angleB=-90,ncurv=1,nodesep=3pt,linearc=.30]{p3}{p4}
$ \end{pspicture}
\end{multline*}
Since by definition
$$\sum_{\pi_r} \kk_{\pi_r}(T^{s_{i_r+1}}, \dots,
T^{s_{i_{r+1}-2}})=\E (T^{s_{i_r+1}}\cdots T^{s_{i_{r+1}-2}})$$ for any
$1\leq r\leq k$, the second part of the proposition follows. \end{proof}

\begin{proof}[Proof of Theorem \ref{theo:formulanapochodna}]
%
Let us consider the case $n=1$. We apply the Leibnitz rule to the
right--hand side of Eq.\ (\ref{eq:rekursjanacyrkularne}):
\begin{multline*}
\frac{d}{dx} \E(T^{s_1} \cdots T^{s_{2m}}) =
\sum_{k\geq 1}
\sum_{\substack{1=i_1<i_2<\cdots\\ \cdots<i_{k+1}=2m+1}}
 \sum_{1\leq r'\leq k} \\
\shoveleft{
\left( \frac{d}{dx} \kk\big( T^{s_{i_{r'}}}, \E(T^{s_{i_{r'}+1}}
T^{s_{i_{r'}+2}} \cdots  T^{s_{i_{r'+1}-2}}) T^{s_{i_{r'+1}-1}} \big) \right)
\times} \\ \times
\prod_{r\neq r'}
\kk\big( T^{s_{i_{r}}}, \E(T^{s_{i_{r}+1}} T^{s_{i_{r}+2}} \cdots
 T^{s_{i_{r+1}-2}}) T^{s_{i_{r+1}-1}} \big).
 \end{multline*}

We denote $j_1=i_{r'}$, $j_2=i_{r'+1}-1$, $k'=r'-1$, $k''=k-r'$, $i'_s=i_s$ for
$1\leq s\leq k'+1$ and $i''_s=i_{s+r'}$ for $1\leq s\leq k''+1$.
Hence
\begin{multline}
\label{eq:wackywacky}
\frac{d}{dx} \E(T^{s_1} \cdots T^{s_{2m}}) =
\sum_{1\leq j_1<j_2\leq 2m} \\
\bigg(  \sum_{k'\geq 1}
\sum_{\substack{1=i'_1<\cdots\\ \cdots<i'_{k'+1}=j_1+1}}
 \prod_{1\leq r \leq k'}
\kk\big( T^{s_{i'_{r}}}, \E(T^{s_{i'_{r}+1}} \cdots
 T^{s_{i'_{r+1}-2}}) T^{s_{i'_{r+1}-1}} \big) \bigg) \times \\
\times \bigg( \frac{d}{dx} \kk\big( T^{s_{j_1}}, \E(T^{s_{j_1+1}}
T^{s_{j_2+2}} \cdots  T^{s_{j_2-1}}) T^{s_{j_2}} \big) \bigg) \times \\
\times \bigg(  \sum_{k''\geq 1}
\sum_{\substack{j_2+1=i''_1<\cdots\\ \cdots<i''_{k''+1}=2m+1}}
 \prod_{1\leq r \leq k''}
 \kk\big( T^{s_{i''_{r}}}, \E(T^{s_{i''_{r}+1}} \cdots
 T^{s_{i''_{r+1}-2}}) T^{s_{i''_{r+1}-1}} \big) \bigg).
\end{multline}
Observe that the first and the third factor on the right--hand side of
\eqref{eq:wackywacky} correspond to the right--hand side of
\eqref{eq:rekursjanacyrkularne}. The second factor can be simplified
by an observation that \eqref{eq:deft} imply $\frac{d}{dx}
\kk(T^{s_{j_1}}, b T^{s_{j_2}})=
b \frac{d}{dx}\kk(T^{s_{j_1}},T^{s_{j_2}})$ for any $b\in\D$. Hence
\begin{multline*}
\frac{d}{dx} \E(T^{s_1} \cdots T^{s_{2m}}) =
\sum_{1\leq j_1<j_2\leq 2m}
\left(\frac{d}{dx} \E(T^{s_{j_1}} T^{s_{j_2}}) \right) \times \\
\times \E(T^{s_1} \cdots T^{s_{j_1-1}})
\E(T^{s_{j_1+1}} \cdots T^{s_{j_2-1}})
\E(T^{s_{j_2+1}} \cdots T^{s_{2m}})
\end{multline*}
which finishes the proof of the case $n=1$.

For the general case we can use recursively \eqref{eq:formulanapochodna} for
$n=1$ and Leibnitz rule in order to compute
$\frac{d^n}{dx^n} \E(T^{s_1} \cdots T^{s_{2m}})$. For example, for
$n=2$ we obtain three summands; for brevity we present below only one of
them:
\begin{multline*}
\frac{d^2}{dx^2} \E(T^{s_1} \cdots T^{s_{2m}})=
\sum_{1\leq j'_1<j'_2\leq 2m}
\left(\frac{d}{dx} \E(T^{s_{j'_1}} T^{s_{j'_2}}) \right) \times \\
\times
\bigg( \sum_{1\leq j''_1<j''_2<j_1'}
\left(\frac{d}{dx} \E(T^{s_{j''_1}} T^{s_{j''_2}}) \right)
\E(T^{s_1} \cdots T^{s_{j''_1-1}})
\E(T^{s_{j''_1+1}} \cdots T^{s_{j''_2-1}}) \times \\ \times
\E(T^{s_{j''_2+1}} \cdots T^{s_{j'_1-1}}) \bigg)
\E(T^{s_{j'_1+1}} \cdots T^{s_{j'_2-1}})
\E(T^{s_{j'_2+1}} \cdots T^{s_{2m}})+ \dots
\end{multline*}
and by renaming the indices we get
\begin{multline*}
\frac{d^2}{dx^2} \E(T^{s_1} \cdots T^{s_{2m}})=
\sum_{0=j_0<\cdots<j_5=2m+1} \\
\shoveleft{
\left( 2 \frac{d}{dx} \E(T^{s_{j_1}} T^{s_{j_2}})
\frac{d}{dx} \E(T^{s_{j_3}} T^{s_{j_4}})+
\frac{d}{dx} \E(T^{s_{j_1}} T^{s_{j_4}})
\frac{d}{dx} \E(T^{s_{j_2}} T^{s_{j_3}}) \right) \times} \\
\times
\prod_{0\leq r\leq 4} \E(T^{s_{j_r+1}} T^{s_{j_r+2}} \cdots
T^{s_{j_{r+1}-1}}).
\end{multline*}
In the general case we get
\begin{multline}
\label{eq:szczegolny}
\frac{d^n}{dx^n} \E(T^{s_1} \cdots T^{s_{2m}})=\\
\sum_{0=j_0<\cdots<j_{2n+1}=2m+1}
c_{s_{j_1},s_{j_2},\dots,s_{j_{2n}}}
\prod_{0\leq r\leq 2n} \E(T^{s_{j_r+1}} T^{s_{j_r+2}} \cdots
T^{s_{j_{r+1}-1}})
\end{multline}
for some constant $c_{s_{j_1},s_{j_2},\dots,s_{j_{2n}}}$ which does not depend
on $m$. This constant can be easily evaluated by the observation that
\eqref{eq:szczegolny} gives in particular
$$\frac{d^n}{dx^n} \E(T^{s_{j_1}} T^{s_{j_2}} \cdots T^{s_{j_{2n}}})=
c_{s_{j_1},s_{j_2},\dots,s_{j_{2n}}}$$
which finishes the proof.

The final remark about nonvanishing summands follows easily from the observation that
if $s_1,\dots,s_{n}\in\{1,-1\}$ and $s_1+\dots+s_n\neq 0$ then every pair partition
of the tuple $(T^{s_1},\dots,T^{s_n})$ must connect $T$ with $T$ or $T\gwia$ with $T\gwia$.
The defining relation  \eqref{eq:deft} implies that $\E(T^{s_1}\cdots T^{s_n})=0$.

\end{proof}

\section{Acknowledgements}
I thank Ken Dykema for many fruitful discussions and for
introducing me into the subject.
Research supported by State Committee for Scientific Research (KBN)
grant No.\ 2 P03A 007 23.
The research was conducted at Texas A\&{}M University on a
scholarship funded by Polish--US Fulbright Commission.
The author is a holder of a scholarship funded by
Fundacja na rzecz Nauki Polskiej.

\bibliography{biblio}

\def\lfhook#1{\setbox0=\hbox{#1}{\ooalign{\hidewidth
  \lower1.5ex\hbox{'}\hidewidth\crcr\unhbox0}}}
\begin{thebibliography}{VDN92}

\bibitem[Bia98]{Biane1998}
Philippe Biane.
\newblock Representations of symmetric groups and free probability.
\newblock {\em Adv. Math.}, 138(1):126--181, 1998.

\bibitem[BLS96]{BozejkoLeinertSpeicher}
Marek Bo{\.z}ejko, Michael Leinert, and Roland Speicher.
\newblock Convolution and limit theorems for conditionally free random
  variables.
\newblock {\em Pacific J. Math.}, 175(2):357--388, 1996.

\bibitem[DH01]{DykemaHaagerup2001}
Kenneth Dykema and Uffe Haagerup.
\newblock Decomposability of {V}oiculescu's circular operator and
  {$DT$}--operators.
\newblock {\em preprint}, 2001.

\bibitem[DY01]{DykemaYan}
Kenneth Dykema and Catherine Yan.
\newblock Generating functions for moments of the quasi--nilpotent {$DT$}
  operator.
\newblock preprint, 2001.

\bibitem[Haa01]{Haagerup2001}
Uffe Haagerup.
\newblock Spectral decomposition of all operators in a {$II_1$} factor, which
  is embedable in {$R^\omega$}.
\newblock preprint, 2001.

\bibitem[KR97a]{KadisonRingrose1}
Richard~V. Kadison and John~R. Ringrose.
\newblock {\em Fundamentals of the theory of operator algebras. {V}ol. {I},
  Reprint of the 1983 original}.
\newblock American Mathematical Society, Providence, RI, 1997.

\bibitem[KR97b]{KadisonRingrose2}
Richard~V. Kadison and John~R. Ringrose.
\newblock {\em Fundamentals of the theory of operator algebras. {V}ol. {I}{I}}.
\newblock American Mathematical Society, Providence, RI, 1997.
\newblock Advanced theory, Corrected reprint of the 1986 original.

\bibitem[Kre72]{Kreweras}
G.~Kreweras.
\newblock Sur les partitions non crois\'ees d'un cycle.
\newblock {\em Discrete Math.}, 1(4):333--350, 1972.

\bibitem[Ran60]{Raney1960}
George~N. Raney.
\newblock Functional composition patterns and power series reversion.
\newblock {\em Trans. Amer. Math. Soc.}, 94:441--451, 1960.

\bibitem[Shl96]{Shlyakhtenko1996}
Dimitri Shlyakhtenko.
\newblock Random {G}aussian band matrices and freeness with amalgamation.
\newblock {\em Internat. Math. Res. Notices}, (20):1013--1025, 1996.

\bibitem[Shl98]{Shlyakhtenko1997}
Dimitri Shlyakhtenko.
\newblock Gaussian random band matrices and operator-valued free probability
  theory.
\newblock In {\em Quantum probability (Gda\'nsk, 1997)}, pages 359--368. Polish
  Acad. Sci., Warsaw, 1998.

\bibitem[SN99]{SpeicherNeu}
Roland Speicher and Peter Neu.
\newblock Physical applications of freeness.
\newblock In {\em XIIth International Congress of Mathematical Physics (ICMP
  '97) (Brisbane)}, pages 261--266. Internat. Press, Cambridge, MA, 1999.

\bibitem[Spe97]{Speicher1997}
Roland Speicher.
\newblock Free probability theory and non-crossing partitions.
\newblock {\em S\'em. Lothar. Combin.}, 39:Art.\ B39c, 38 pp.\ (electronic),
  1997.

\bibitem[Spe98]{Speicher1998}
Roland Speicher.
\newblock Combinatorial theory of the free product with amalgamation and
  operator-valued free probability theory.
\newblock {\em Mem. Amer. Math. Soc.}, 132(627):x+88, 1998.

\bibitem[VDN92]{VoiculescuDykemaNica}
D.~V. Voiculescu, K.~J. Dykema, and A.~Nica.
\newblock {\em Free random variables}.
\newblock American Mathematical Society, Providence, RI, 1992.
\newblock A noncommutative probability approach to free products with
  applications to random matrices, operator algebras and harmonic analysis on
  free groups.

\bibitem[Voi91]{Voiculescu1991}
Dan Voiculescu.
\newblock Limit laws for random matrices and free products.
\newblock {\em Invent. Math.}, 104(1):201--220, 1991.

\bibitem[Voi95]{VoiculescuPlenar}
Dan Voiculescu.
\newblock Free probability theory: random matrices and von {N}eumann algebras.
\newblock In {\em Proceedings of the International Congress of Mathematicians,
  Vol.\ 1, 2 (Z\"urich, 1994)}, pages 227--241, Basel, 1995. Birkh\"auser.

\bibitem[Voi00]{VoiculescuLectures}
Dan Voiculescu.
\newblock Lectures on free probability theory.
\newblock In {\em Lectures on probability theory and statistics (Saint-Flour,
  1998)}, pages 279--349. Springer, Berlin, 2000.

\end{thebibliography}

\end{document}